\documentclass[fleqn]{mat01}
\usepackage{times,mathtimy,amssymb,latexsym}
\begin{document}

\setcounter{page}{97} \firstpage{97}

\renewcommand\theequation{\thesection\arabic{equation}}

\def\theorre{\trivlist \item[\hskip \labelsep{\bf Theorem MYX \cite{19}.}]}
\def\rema{\trivlist \item[\hskip \labelsep{\it Remark.}]}

\newtheorem{theo}[defin]{\bf Theorem}
\newtheorem{lem}[defin]{Lemma}
\newtheorem{rem}[defin]{\it Remark}

\markboth{Yongfu Su and Xiaolong Qin}{Strong convergence for
nonlinear mappings}

\title{Strong convergence of modified Ishikawa
iterations\\ for nonlinear mappings}

\author{YONGFU SU and XIAOLONG QIN$^{*}$}

\address{Department of Mathematics, Tianjin Polytechnic University, Tianjin 300160,
China\\
\noindent $^{*}$E-mail: qxlxajh@163.com}

\volume{117}

\mon{February}

\parts{1}

\pubyear{2007}

\Date{MS received 25 May 2006}

\begin{abstract}
In this paper, we prove a strong convergence theorem of modified
Ishikawa iterations for relatively asymptotically nonexpansive
mappings in Banach space. Our results  extend and improve the
recent results by Nakajo, Takahashi, Kim, Xu, Matsushita and some
others.
\end{abstract}


\maketitle

\section{Introduction and preliminaries}

Let $E$ be a real Banach space, $C$ a nonempty closed convex
subset of $E$ and $T\hbox{\rm :}\ C\rightarrow C$ a mapping.
Recall that $T$ is nonexpansive if
\begin{equation*}
\|Tx-Ty\|\leq \|x-y\| \quad \mbox{for all}\ x,y\in C,
\end{equation*}
and $T$ is asymptotically nonexpansive \cite{10} if there exists a
sequence $\{k_n\}$ of positive real numbers with
$\lim_{n\rightarrow\infty}k_n=1$ such that
\begin{equation*}
\|T^nx-T^ny\|\leq k_n\|x-y\| \quad\mbox{for all}\ n\geq1\
\mbox{and}\ x,y\in C.
 \end{equation*}
A point $x\in C$ is a fixed point of $T$ provided $Tx=x$. Denote
by $F(T)$ the set of fixed points of $T$; that is, $F(T)=\{x\in
C\hbox{\rm :}\ Tx=x\}$.

Some iteration processes are often used to approximate a fixed
point of a nonexpansive mapping.  The first iteration process is
now known as Mann's iteration process \cite{14} which is defined
as
\begin{equation}
x_{n+1}=\alpha_nx_n+(1-\alpha_n)Tx_n,\ \ \ n\geq 0
\end{equation}
where the initial guess $x_0$ is taken in $C$ arbitrarily and the
sequence $\{\alpha_n\}_{n=0}^\infty$ is in the interval [0, 1].

The second iteration process is referred to as Ishikawa's
iteration process \cite{11} which is defined recursively by
\begin{equation}
\begin{cases}
y_n = \beta_nx_n+(1-\beta_n)Tx_n, \\
x_{n+1} = \alpha_nx_n+(1-\alpha_n)Ty_n,
\end{cases}
\end{equation}
where the initial guess $x_0$ is taken in $C$ arbitrarily and
$\{\alpha_n\}$ and $\{\beta_n\}$ are sequences in the interval [0,
1].

In general, not much has been known about the convergence of the
iteration processes (1.1) and (1.2) unless the underlying space
$E$ has elegant properties which we briefly mention here.

Reich \cite{18} proved that if $E$ is a uniformly convex Banach
space with a Fr\'{e}chet differentiable norm and if $\{\alpha_n\}$
is chosen such that
$\sum_{n=1}^\infty\alpha_n(1-\alpha_n)=\infty$, then the sequence
$\{x_n\}$ defined by (1.1) converges weakly to a fixed point of
$T$. However we note that Mann's iterations have only weak
convergence even in a Hilbert space \cite{9}.

Attempts to modify the Mann's iteration method (1.1) so that
strong convergence is guaranteed have recently been made. Nakajo
and Takahashi \cite{15} proposed the following modification of the
Mann's iteration (1.1) for a single nonexpansive mapping $T$ in a
Hilbert space:
\begin{equation}
\begin{cases}
x_0\in C \ \hbox{arbitrarily\ chosen},\\[.2pc]
y_n = \alpha_nx_n+(1-\alpha_n)Tx_n, \\[.2pc]
C_n=\{z\in C\hbox{\rm :}\ \|y_n-z\|\leq \|x_n-z\|\},\\[.2pc]
Q_n=\{z\in C\hbox{\rm :}\ \langle x_0-x_n, x_n-z\rangle\geq 0\},\\[.2pc]
x_{n+1}=P_{C_n\cap Q_n}x_0,
\end{cases}
\end{equation}
where $P_K$ denotes the metric projection from $H$ onto a closed
convex subset $K$ of $H$ and proved that sequence $\{x_n\}$
converges strongly to $P_{F(T)}x_0$.

Recently, Kim and Xu \cite{13} has adapted the iteration (1.1) in
a Hilbert space. More precisely, they introduced the following
iteration process for asymptotically nonexpansive mappings, with
$C$ a closed convex bounded subset of a Hilbert space:
\begin{equation}
\begin{cases}
x_0\in C \ \hbox{arbitrarily\ chosen,}\\[.2pc]
y_n = \alpha_nx_n+(1-\alpha_n)T^nx_n, \\[.2pc]
C_n=\{z\in C\hbox{\rm :}\ \|y_n-z\|^2\leq \|x_n-z\|^2+\theta_n\},\\[.2pc]
Q_n=\{z\in C\hbox{\rm :}\ \langle x_0-x_n, x_n-z\rangle\geq 0\},\\[.2pc]
x_{n+1}=P_{C_n \cap Q_n}x_0,
\end{cases}
\end{equation}
where
\begin{equation*}
\theta_n=(1-\alpha_n)(k_n^2-1)(\hbox{diam}\ C)^2\rightarrow 0\ \ \
\mbox{as}\ n\rightarrow\infty.
\end{equation*}
They proved $\{x_n\}$ converges in norm to $P_{F(T)}x_0$ under the
conditions that the sequence $\{\alpha_n\}_{n=0}^\infty$ in $(0,
1)$ such that $\alpha_n\leq \alpha$ for all $n$ and for some
$0<\alpha<1.$

On the other hand, process (1.2) is indeed more general than
process (1.1). But research has been done on the latter due to
reasons that the formulation of process (1.1) is simpler than that
of (1.2) and that a convergence theorem for process (1.1) may lead
to a convergence theorem for process (1.2) provided that
$\{\beta_n\}$ satisfies certain appropriate conditions. However,
the introduction of the process (1.2) has its own right. Actually,
the process (1.1) may fail to converge while process (1.2) can
still converge for a Lipschitz pseudo-contractive mapping in a
Hilbert space \cite{7}.

In \cite{19}, Martinez-Yanes and Xu proved the following theorem.

\begin{theorre}
{\it Let $C$ be a closed convex subset of a Hilbert space H and
let $T\hbox{\rm :}\ \rightarrow C$ be a nonexpansive mapping such
that $F(T)\neq \emptyset$. Assume that $\{\alpha_n\}_{n=0}^\infty$
and $\{\beta_n\}_{n=0}^\infty$ are sequences in $[0, 1]$ such that
$\alpha_n\leq 1-\delta$ for some $\delta\in (0, 1]$ and
$\beta_n\rightarrow 1$. Define a sequence $\{x_n\}_{n=0}^\infty$
in $C$ by the algorithm{\rm :}
\begin{equation*}
\begin{cases}
x_0\in C \quad\mbox{arbitrarily\ chosen}\\[.2pc]
z_n=\beta_nx_n+(1-\beta_n)Tx_n,\\[.2pc]
y_n =\alpha_nx_n+(1-\alpha_n)Tz_n, \\[.2pc]
C_n=\{v\in C\hbox{\rm :}\ \|y_n-v\|^2\leq\|x_n-v\|^2\\[.2pc]
\qquad +(1-\alpha_n)(\|z_n\|^2-\|x_n\|^2+2\langle x_n-z_n, v\rangle)\},\\[.2pc]
Q_n=\{v\in C\hbox{\rm :}\ \langle x_0-x_n, x_n-v\rangle\geq 0\},\\[.2pc]
x_{n+1}=\Pi_{C_n \cap Q_n}x_0.
\end{cases}
\end{equation*}
Then $\{x_n\}$ converges in norm to $P_{F(T)}x_0.$}\vspace{.5pc}
\end{theorre}

The purpose of this paper is to employ Nakajo and Takahashi's idea
\cite{15} to modify process (1.2) for relatively asymptotically
nonexpansive mappings to have strong convergence theorem in Banach
spaces.

Let $E$ be a smooth Banach space with dual $E^*$. We denote by $J$
the normalized duality mapping from $E$ to $2^{E^*}$ defined by
\begin{equation*}
Jx=\{f^*\in E^*\hbox{\rm :}\ \langle
x,f^*\rangle=\|x\|^2=\|f^*\|^2\},
\end{equation*}
where $\langle \cdot,\cdot\rangle$ denotes the generalized duality
pairing.

As we all know that if $C$ is a nonempty closed convex subset of a
Hilbert space $H$ and $P_C\hbox{\rm :}\ H\rightarrow C$ is the
metric projection of $H$ onto $C$, then $P_C$ is nonexpansive.
This fact actually characterizes Hilbert spaces and consequently,
it is not available in more general Banach spaces. In this
connection, Alber [1] recently introduced a generalized projection
operator $\Pi_C$ in a Banach space $E$ which is an analogue of the
metric projection in Hilbert spaces.

Consider the functional defined by
\begin{equation}
\phi(x,y)=\|x\|^2-2\langle x,j(y)\rangle+\|y\|^2\quad \mbox{for}\
\ x,y\in E,
\end{equation}
where $j(y)\in J(y).$ Observe that, in a Hilbert space
$H$, (1.5) reduces to $\phi(x,y)=\|x-y\|^2,\ \ x,y\in H.$

The generalized projection $\Pi_C\hbox{\rm :}\ E\rightarrow C$ is
a map that assigns to an arbitrary point $x\in E$ the minimum
point of the functional $\phi(x,y),$ that is, $\Pi_Cx=\bar{x},$
where $\bar{x}$ is the solution to the minimization problem
\begin{equation}
\phi(\bar{x},x)=\inf\limits_{y\in C}\phi(y,x).
\end{equation}
The existence and uniqueness of the operator $\Pi_C$ follow from
the properties of the functional $\phi(x,y)$ and strict
monotonicity of the mapping $J$ (see, for example, \cite{3}). In
Hilbert spaces, $\Pi_C=P_C.$ It is obvious from the definition of
function $\phi$ that
\begin{equation}
(\|y\|-\|x\|)^2\leq\phi(y,x)\leq(\|y\|+\|x\|)^2\quad \mbox{for
all}\ x,y\in E.
\end{equation}

\begin{rema}
If $E$ is a strictly convex and smooth Banach space, then for
$x,y\in E$, $\phi(x,y)=0$ if and only if $x=y$. It is sufficient
to show that if $\phi(x,y)=0$ then $x=y$. From (1.7), we have
$\|x\|=\|y\|$. This implies $\langle
x,Jy\rangle=\|x\|^2=\|Jy\|^2.$ From the definitions of $J,$ we
have $Jx=Jy$. Since $J$ is one-to-one, we have $x=y;$ see [8,19]
for more details.

Let $C$ be a closed convex subset of $E$, and let $T$ be a mapping
from $C$ into itself. We denote by $F(T)$ the set of fixed points
of $T$. A point of $p$ in $C$ is said to be an asymptotically
fixed point of $T$ [17] if $C$ contains a sequence $\{x_n\}$ which
converges weakly to $p$ such that the strong
$\lim_{n\rightarrow\infty}(Tx_n-x_n)=0.$ The set of asymptotic
fixed points of $T$ will be denoted by $\hat{F}(T)$. A mapping $T$
from $C$ into itself is called relatively nonexpansive
\cite{1,2,3} if $\hat{F}(T)=F(T)$ and $\phi(p,Tx)\leq \phi(p,x)$
for all $x\in C$ and $p\in F(T)$. A mapping $T$ from $C$ into
itself is called relatively asymptotically nonexpansive if
$\hat{F}(T)=F(T)$ and $\phi(p,T^nx)\leq k_n^2\phi(p,x)$ for all
$x\in C$ and $p\in F(T)$.

A Banach space $E$ is said to be strictly convex if
$\|\frac{x+y}{2}\|<1$ for all $x,y\in E$ with $\|x\|=\|y\|=1$ and
$x\neq y$. It is said to be uniformly convex if
$\lim_{n\rightarrow\infty}\|x_n-y_n\|=0$ for any two sequences
$\{x_n\},\ \{y_n\}$ in $E$ such that $\|x_n\|=\|y_n\|=1$ and
$\lim_{n\rightarrow\infty}\|\frac{x_n+y_n}{2}\|=1$. Let $U=\{x\in
E\hbox{\rm :}\ \|x\|=1\}$ be the unit sphere of $E$. Then the
Banach space $E$ is said to be smooth provided
\begin{equation*}
\lim\limits_{t\rightarrow 0}\frac{\|x+ty\|-\|x\|}{t}
\end{equation*}
exists for each $x,y\in U.$ It is also said to be uniformly smooth
if the limit is attained uniformly for $x,y\in E$. It is
well-known that if $E$ is uniformly smooth, then $J$ is uniformly
norm-to-norm continuous on each bounded subset of $E$. A Banach
space is said to have the Kadec--Klee property if a sequence
$\{x_n\}\rightharpoonup x\in E$ and $\|x_n\|\rightarrow \|x\|$,
then $x_n\rightarrow x$. It is known that if $E$ is uniformly
convex then $E$ has the Kadec--Klee property; see \cite{10,19} for
more details.
\end{rema}

We need the following Lemmas for the proof of our main results.

\begin{lem}\hskip -.5pc {\rm \cite{12}.} \ \ Let $E$ be a uniformly convex and smooth Banach space and let
$\{x_n\},$ $\{y_n\}$ be two sequences of $E$. If
$\phi(x_n,y_n)\rightarrow 0$ and either $\{x_n\}$ or $\{y_n\}$ is
bounded{\rm ,} then $x_n-y_n\rightarrow 0.$
\end{lem}

\begin{lem}\hskip -.5pc {\rm \cite{1,2,3}.} \ \ Let $C$ be a nonempty closed convex subset of a smooth
Banach space $E$ and $x\in E$. Then{\rm ,} $x_0=\Pi_Cx$ if and
only if
\begin{equation*}
\langle x_0-y,Jx-Jx_0\rangle\geq 0\quad \mbox{for}\ y\in C.
\end{equation*}
\end{lem}

\begin{lem}\hskip -.5pc {\rm \cite{1,2,3}.} \ \ Let $E$ be a reflexive{\rm ,} strictly
convex and smooth Banach space. Let $C$ be a nonempty closed
convex subset of $E$ and let $x\in E.$ Then
\begin{equation*}
\phi(y,\Pi_cx)+\phi(\Pi_cx,x)\leq\phi(y,x)\quad \mbox{for all}\
y\in C.
\end{equation*}
\end{lem}

\begin{lem}
Let $E$ be a uniformly convex and uniformly smooth Banach space.
Let $C$ be a closed convex subset of $E$ and let $T$ be a
relatively asymptotically nonexpansive mapping from $C$ into
itself. If $T$ is continuous, then $F(T)$ is closed and convex.
\end{lem}

\begin{proof}
We first show that $F(T)$ is closed. Since  $T$ is continuous, we
can obtain the closedness of $F(T)$ easily. Next, we show that
$F(T)$ is convex for $x,y\in F(T)$ and $t\in (0,1)$. Put
$p=tx+(1-t)y.$ It is sufficient to show $Tp=p.$ In fact, we
have\pagebreak
\begin{align*}
\phi(p,T^np)&=\|p\|^2-2\langle p,JT^np\rangle+\|T^np\|^2\\[.4pc]
&=\|p\|^2-2\langle tx+(1-t)y,JT^np\rangle+\|T^np\|^2\\[.4pc]
&=\|p\|^2-2t\langle x,JT^np\rangle-2(1-t)\langle y,
JT^np\rangle+\|T^np\|^2\\[.4pc]
&=\|p\|^2+t\phi(x, T^np)+(1-t)\phi(y,T^np)-t\|x\|^2\\
&\quad\,-(1-t)\|y\|^2\\[.4pc]
&\leq\|p\|^2+k_nt\phi(x, p)+k_n(1-t)\phi(y,p)-t\|x\|^2\\
&\quad\,-(1-t)\|y\|^2\\[.4pc]
&=(k_n-1)(t\|x\|^2+(1-t)\|y\|^2-\|p\|^2).
\end{align*}
Take the limit as $n\rightarrow\infty$ yields
\begin{equation*}
\lim_{n\rightarrow\infty}\phi(p,T^np)=0.
\end{equation*}
Now we apply Lemma 1.1 to see that $T^np\rightarrow p$ strongly.
By continuity of $T$ we obtain $p\in F(T).$ This completes the
lemma 1.4.\hfill $\Box$
\end{proof}

\section{Main results}

\setcounter{equation}{0}

\setcounter{defin}{0}
\begin{theo}[\!]
Let $E$ be a uniformly convex and uniformly smooth Banach space.
Let $C$ be a nonempty bounded closed convex subset of $E$. Let
$T\hbox{\rm :}\ C\rightarrow C$ be a relatively asymptotically
nonexpansive mapping with sequence $\{k_n\}$ such that
$k_n\rightarrow 1$ as $n\rightarrow\infty$ and
$F(T)\neq\emptyset$. Assume that $\{\alpha_n\}_{n=0}^\infty$ and
$\{\beta_n\}_{n=0}^\infty$ are sequences in $[0,1]$ such that
$\limsup_{n\rightarrow\infty}\alpha_n<1$ and $\beta_n\rightarrow
1.$ Define a sequence $\{x_n\}$ in $C$ by the following
algorithm{\rm :}
\begin{equation}
\begin{cases}
x_0\in C \ \ \ {\rm arbitrarily\ chosen,}\\[.2pc]
z_n=J^{-1}(\beta_nJx_n+(1-\beta_n)JT^nx_n),\\[.2pc]
y_n = J^{-1}(\alpha_nJx_n+(1-\alpha_n)JT^nz_n), \\[.2pc]
C_n=\{v\in C\hbox{\rm :}\ \phi(v,y_n)\leq\phi(v,x_n)\\
\qquad+(1-\alpha_n)(k_n^2\|z_n\|^2-\|x_n\|^2+(k_n^2-1)M-2\langle v, k_n^2Jz_n-Jx_n\rangle)\},\\[.2pc]
Q_n=\{v\in C\hbox{\rm :}\ \langle Jx_0-Jx_n, x_n-v\rangle\geq 0\},\\[.2pc]
x_{n+1}=\Pi_{C_n \cap Q_n}x_0,
\end{cases}
\end{equation}
where $J$ is the duality mapping on $E$ and  $M$ is an appropriate
constant such that $M>\|v\|^2$ for each $v\in C.$ If $T$ is
uniformly  continuous{\rm ,} then $\{x_n\}$ converges to some $q=
P_{F(T)}x_0.$
\end{theo}

\begin{proof}
We first show that $C_n$ and $Q_n$ are closed and convex for each
$n\geq 0$. From the definition of $C_n$ and $Q_n$, it is obvious
that $C_n$ is closed and $Q_n$ is closed and convex for each
$n\geq 0.$ We prove that $C_n$ is convex. For $v_1, v_2\in C_n$
and $t\in (0,1),$ put $v=tv_1+(1-t)v_2$. It is sufficient to show
that $v\in C_n.$ Since
\begin{align*}
\phi(v,y_n)&\leq\phi(v,x_n)+(1-\alpha_n)(k_n^2\|z_n\|^2-\|x_n\|^2\\[.2pc]
&\quad\,+(k_n^2-1)M-2\langle v, k_n^2Jz_n-Jx_n\rangle)
\end{align*}
is equivalent to
\begin{align*}
&2\langle v, Jx_n\rangle+2(1-\alpha_n)\langle
v,k_n^2Jz_n-Jx_n\rangle-2\langle v,Jy_n\rangle\\[.2pc]
&\leq(2-\alpha_n)\|x_n\|^2+(1-\alpha_n)(k_n^2\|z_n\|^2+(k_n^2-1)M)-\|y_n\|^2,
\end{align*}
we have
\begin{align*}
 &2\langle v,
Jx_n\rangle+2(1-\alpha_n)\langle v,k_n^2Jz_n-Jx_n\rangle-2\langle v,Jy_n\rangle\\
&\quad\,=2\langle tv_1+(1-t)v_2,
Jx_n\rangle+2(1-\alpha_n)\langle tv_1+(1-t)v_2,k_n^2Jz_n-Jx_n\rangle\\
&\qquad\, -2\langle tv_1+(1-t)v_2,Jy_n\rangle\\[.4pc]
&\quad\,=2t\langle v_1, Jx_n\rangle+2(1-t)\langle v_2,
Jx_n\rangle+2(1-\alpha_n)t\langle v_1,k_n^2Jz_n-Jx_n\rangle\\
&\qquad\, +2(1-\alpha_n)(1-t)\langle
v_2,k_n^2Jz_n-Jx_n\rangle\\
&\qquad\,-2t\langle v_1,Jy_n\rangle-2(1-t)\langle
v_2,Jy_n\rangle\\[.4pc]
&\quad\,
\leq(2-\alpha_n)\|x_n\|^2+(1-\alpha_n)(k_n^2\|z_n\|^2+(k_n^2-1)M)-\|y_n\|^2.
\end{align*}
This implies $v\in C_n$. So $C_n$ is convex. Next, we show that
$F(T)\subset C_n$ for all $n$. Indeed, we have, for all $p\in
F(T)$,
\begin{align*}
\phi(p,y_n)&=\phi(p, J^{-1}(\alpha_nJx_n+(1-\alpha_n)JT^nz_n))\\[.4pc]
&=\|p\|^2-2\langle
p,\alpha_nJx_n+(1-\alpha_n)JT^nz_n\rangle\\
&\quad\,+\|\alpha_nJx_n+(1-\alpha_n)JT^nz_n)\|^2\\[.4pc]
&\leq\|p\|^2-2\alpha_n\langle p,Jx_n\rangle-2(1-\alpha_n)\langle
p,JT^nz_n\rangle\\
&\quad\,
+\alpha_n\|x_n\|^2+(1-\alpha_n)\|T^nz_n\|^2\\[.4pc]
&\leq\alpha_n\phi(p,x_n)+(1-\alpha_n)\phi(p,T^nz_n)\\[.4pc]
&\leq\alpha_n\phi(p,x_n)+k_n^2(1-\alpha_n)\phi(p,z_n)\\[.4pc]
&=\phi(p,x_n)+(1-\alpha_n)[k_n^2\phi(p,z_n)-\phi(p,x_n)]\\[.4pc]
&\leq\phi(p,x_n)+(1-\alpha_n)(k_n^2\|z_n\|^2-\|x_n\|^2\\
&\quad\,+(k_n^2-1)\|p\|^2-2\langle p, k_n^2Jz_n-Jx_n\rangle)\\[.4pc]
&\leq\phi(p,x_n)+(1-\alpha_n)(k_n^2\|z_n\|^2-\|x_n\|^2\\
&\quad\,+(k_n^2-1)M-2\langle p, k_n^2Jz_n-Jx_n\rangle).
\end{align*}
So $p\in C_n$ for all $n.$ Next we show that
\begin{equation}
F(T)\subset Q_n,\quad \mbox{for all}\ n\geq 0.
\end{equation}
We prove this by induction. For $n=0,$ we have $F(T)\subset
C=Q_0.$ Assume that $F(T)\subset Q_n.$ Since $x_{n+1}$ is the
projection of $x_0$ onto $C_n\cap Q_n$, by Lemma 1.2 we have
\begin{equation*}
\langle Jx_0-Jx_{n+1}, x_{n+1}-z\rangle\geq 0,\quad \forall z\in
C_n\cap Q_n.
\end{equation*}
As $F(T)\subset C_n\cap Q_n$ by the induction assumptions, the
last inequality holds, in particular, for all $z\in F(T)$. This
together with the definition of $Q_{n+1}$ implies that
$F(T)\subset Q_{n+1}.$ Hence (2.2) holds for all $n\geq 0.$ This
implies that $\{x_n\}$  is well defined. Since
$x_{n+1}=\Pi_{C_n\cap Q_n}x_0\in Q_n$, we have
\begin{equation*}
\phi(x_n,x_0)\leq \phi(x_{n+1},x_0)\quad \mbox{for all}\ n\geq 0.
\end{equation*}
Therefore $\{\phi(x_n,x_0)\}$ is nondecreasing. Since $C$ is
bounded, $\phi(x_n,x_0)$ is bounded. Moreover from (1.7), we have
that $\{x_n\}$ is bounded. So, we obtain that the limit of
$\{\phi(x_n,x_0)\}$ exists. From Lemma 1.3, we have
\begin{align*}
\phi(x_{n+1},x_n)&=\phi(x_{n+1},\Pi_{C_n}x_0)\leq\phi(x_{n+1},x_0)-\phi(\Pi_{C_n}x_0,x_0)\\[.4pc]
&=\phi(x_{n+1},x_0)-\phi(x_n,x_0)
\end{align*}
for all $n\geq 0.$ This implies that
\begin{equation}
\lim_{n\rightarrow\infty}\phi(x_{n+1},x_n)=0.
\end{equation}
By using Lemma 1.1, we obtain
\begin{equation}
\lim\limits_{n\rightarrow\infty}\|x_{n+1}-x_n\| =0.
\end{equation}
Since $x_{n+1}=\Pi_{C_n\cap Q_n}x_0\in C_n$, from the definition
of $C_n,$ we also have
\begin{align}
\phi(x_{n+1},y_n) &\leq
\phi(x_{n+1},x_n)+(1-\alpha_n)(k_n^2\|z_n\|^2-\|x_n\|^2\nonumber\\[.2pc]
&\quad\,+(k_n^2-1)M-2\langle x_{n+1}, k_n^2Jz_n-Jx_n\rangle).
\end{align}
However, since $\lim_{n\rightarrow\infty}\beta_n=1$ and $\{x_n\}$
is bounded, we obtain
\begin{align*}
\phi(x_{n+1},z_n)&=\phi(x_{n+1},J^{-1}(\beta_nJx_n+(1-\beta_n)JT^nx_n))\\[.4pc]
&=\|x_{n+1}\|^2-2\langle
x_{n+1},\beta_nJx_n+(1-\beta_n)JT^nx_n)\rangle\\
&\quad\,
+\|\beta_nJx_n+(1-\beta_n)JT^nx_n\|^2\\[.4pc]
&\leq \|x_{n+1}\|^2-2\beta_n\langle
x_{n+1},Jx_n\rangle-2(1-\beta_n)\langle x_{n+1},JT^nx_n\rangle\\
&\quad\,
+\beta_n\|x_n\|^2+(1-\beta_n)\|T^nx_n\|^2 \\[.4pc]
&=\beta_n\phi(x_{n+1},x_n)+(1-\beta_n)\phi(x_{n+1},T^nx_n).
\end{align*}
Therefore, we obtain
\begin{equation}
\phi(x_{n+1},z_n)\rightarrow 0,
\end{equation}
which yields
\begin{equation}
\|x_{n+1}\|^2+\|z_n\|^2-2\langle x_{n+1}, Jz_n\rangle\rightarrow
0. \end{equation}
On the other hand, we have
\begin{align}
&k_n^2\|z_n\|^2-\|x_n\|^2-2\langle x_{n+1},
k_n^2Jz_n-Jx_n\rangle\nonumber\\[.4pc]
&\quad\,= k_n^2\|z_n\|^2-\|x_n\|^2-2k_n^2\langle x_{n+1},
Jz_n\rangle+2\langle x_{n+1}, Jx_n\rangle\nonumber\\[.4pc]
&\quad\, = (k_n^2\|z_n\|^2+k_n^2\|x_{n+1}\|^2-2k_n^2\langle
x_{n+1}, Jz_n\rangle)\nonumber\\
&\qquad\,+2\langle x_{n+1},
Jx_n\rangle-k_n^2\|x_{n+1}\|^2-\|x_n\|^2.
\end{align}
Now, we consider
\begin{align*}
&2\langle x_{n+1}, Jx_n\rangle-k_n^2\|x_{n+1}\|^2-\|x_n\|^2\\[.4pc]
&\quad\,=\langle x_{n+1}, Jx_n\rangle+\langle x_{n+1}, Jx_n\rangle-k_n^2\|x_{n+1}\|^2-\|x_n\|^2\\[.4pc]
&\quad\,=\langle x_n+x_{n+1}-x_n, Jx_n\rangle+\langle x_{n+1}, Jx_{n+1}+Jx_n-Jx_{n+1}\rangle\\
&\qquad\,-k_n^2\|x_{n+1}\|^2-\|x_n\|^2\\[.4pc]
&\quad\,= \langle x_{n+1}-x_n, Jx_n\rangle+\langle x_{n+1},
Jx_n-Jx_{n+1}\rangle-(k_n^2-1)\|x_{n+1}\|^2.
\end{align*}
It follows from (2.4) that
\begin{equation}
2\langle x_{n+1},
Jx_n\rangle-k_n^2\|x_{n+1}\|^2-\|x_n\|^2\rightarrow 0.
\end{equation}
It follows from (2.7) and (2.9) that
\begin{equation}
k_n^2\|z_n\|^2-\|x_n\|^2-2\langle x_{n+1},
k_nJz_n-Jx_n\rangle\rightarrow 0.
\end{equation}
Combining (2.3), (2.5) and (2.10), we have
\begin{equation*}
\lim\limits_{n\rightarrow\infty}\phi(x_{n+1},y_n)=0.
\end{equation*}
Using Lemma 1.1, we obtain
\begin{equation}
\lim\limits_{n\rightarrow\infty}\|x_{n+1}-y_n\| =0.
\end{equation}
Since $J$ is uniformly norm-to-norm continuous on bounded sets, we
have
\begin{equation}
\lim\limits_{n\rightarrow\infty}\|Jx_{n+1}-Jy_n\|
=\lim\limits_{n\rightarrow\infty}\|Jx_{n+1}-Jx_n\|=0.
\end{equation}
Notice that
\begin{align*}
\|Jx_{n+1}-Jy_n\|&=\|Jx_{n+1}-(\alpha_nJx_n+(1-\alpha_n)JT^nz_n)\|\\[.4pc]
&=\|\alpha_n(Jx_{n+1}-Jx_n)+(1-\alpha_n)(Jx_{n+1}-JT^nz_n)\|\\[.4pc]
&=\|(1-\alpha_n)(Jx_{n+1}-JT^nz_n)-\alpha_n(Jx_n-Jx_{n+1})\|\\[.4pc]
&\geq(1-\alpha_n)\|Jx_{n+1}-JT^nz_n\|-\alpha_n\|Jx_n-Jx_{n+1}\|.
\end{align*}
We have
\begin{align*}
\|Jx_{n+1}-JT^nz_n\|\leq\frac{1}{1-\alpha_n}(\|Jx_{n+1}-Jy_n\|+\alpha_n\|Jx_n-Jx_{n+1}\|).
\end{align*}
From (2.12) and $\limsup_{n\rightarrow\infty}\alpha_n<1$, we
obtain
\begin{equation*}
\lim\limits_{n\rightarrow\infty}\|Jx_{n+1}-JT^nz_n\|=0.
\end{equation*}
Since $J^{-1}$ is also uniformly norm-to-norm continuous on
bounded sets, we obtain
\begin{equation}
\lim\limits_{n\rightarrow\infty}\|x_{n+1}-T^nz_n\|=0
\end{equation}
and hence
\begin{equation*}
\|x_n-T^nx_n\|\leq\|x_{n+1}-x_n\|+\|x_{n+1}-T^nz_n\|.
\end{equation*}
It follows from (2.4) and (2.13) that
$\lim_{n\rightarrow\infty}\|T^nx_n-x_n\|=0.$ Putting
$L\!=\!\sup\{k_n\hbox{\rm :}\ n\!\geq\! 1\}$ $<\infty,$ we obtain
\begin{align*}
\|Tx_n-x_n\|&\leq\|Tx_n-T^{n+1}x_n\|+\|T^{n+1}x_n-T^{n+1}x_{n+1}\|\\[.2pc]
&\quad\, +\|T^{n+1}x_{n+1}-x_{n+1}\|+\|x_{n+1}-x_n\|.
\end{align*}
Since $T$ is uniformly  continuous, we have
\begin{equation*}
\|Tx_n-x_n\|\rightarrow 0, \quad \mbox{as}\ n\rightarrow\infty.
\end{equation*}
Finally, we prove that $x_n\rightarrow q= \Pi_{F(T)}x_0$. Assume
that $\{x_{n_i}\}$ is a subsequence of $\{x_n\}$ such that
$\{x_{n_i}\}\rightharpoonup q\in C$. Then $q\in\hat{F}(T)=F(T)$.
Next we show that $q= \Pi_{F(T)}x_0$ and convergence is strong.
Putting $q'= \Pi_{F(T)}x_0$ from $x_{n+1}=\Pi_{C_n\cap Q_n}x_0$
and $q'\in F(T)\subset C_n\cap Q_n$. We have
$\phi(x_{n+1},x_0)\leq\phi(q',x_0)$. On the other hand, from
weakly lower semicontinuity  of the norm, we obtain
\begin{align*}
\phi(q,x_0)&=\|q\|^2-2\langle
q,Jx_0\rangle+\|x_0\|^2\\[.4pc]
&\leq\liminf\limits_{i\rightarrow\infty}(\|x_{n_i}\|^2-\langle
x_{n_i}, Jx_0\rangle+\|x_0\|^2)\\[.4pc]
&\leq\liminf\limits_{i\rightarrow\infty}\phi(x_{n_i},x_0)\leq\limsup\limits_{i\rightarrow\infty}\phi(x_{n_i},x_0)\\[.4pc]
&\leq\phi(q',x_0).
\end{align*}
It follows from the definition of $\Pi_{F(T)}x_0$, that $q=
\Pi_{F(T)}x_0$ and hence
\begin{equation*}
\lim\limits_{i\rightarrow\infty}\phi(x_{n_i},x_0)=\phi(q',x_0)=\phi(q,
x_0).
\end{equation*}
So we have $\lim_{i\rightarrow\infty}\|x_{n_i}\|=\|q\|$. Using the
Kadec--Klee property of $E$, we obtain that $\{x_{n_i}\}$
converges strongly to $q= P_{F(T)}x_0$. Since $\{x_{n_i}\}$ is an
arbitrarily weakly convergent sequence of $\{x_n\}$, we can
conclude that $\{x_n\}$ converges strongly to $\Pi_{F(T)}x_0.$
This completes the proof.\hfill $\Box$
\end{proof}

\section{Applications}

\setcounter{defin}{0}

\begin{theo}[\!]
Let $C$  be a nonempty bounded closed convex subset of a Hilbert
space $H$ and let $T\hbox{\rm :}\ C\rightarrow C$ be an
asymptotically nonexpansive mapping with sequence $\{k_n\}$ such
that $k_n\rightarrow 1$ as $n\rightarrow\infty$. Assume that
$\{\alpha_n\}_{n=0}^\infty$ and $\{\beta_n\}_{n=0}^\infty$ are
sequences in $[0,1]$ such that
$\limsup_{n\rightarrow\infty}\alpha_n<1$ and $\beta_n\rightarrow
1.$ Define a sequence $\{x_n\}$ in $C$ by the following
algorithm{\rm :}
\begin{equation*}
\begin{cases}
x_0\in C \quad \textit{arbitrarily\ chosen},\\
z_n=\beta_nx_n+(1-\beta_n)T^nx_n,\\
y_n = \alpha_nx_n+(1-\alpha_n)T^nz_n, \\
C_n=\{v\in C\hbox{\rm :}\ \|y_n-v\|^2\leq\|x_n-v\|^2\\
\qquad\,+(1-\alpha_n)(k_n^2\|z_n\|^2-\|x_n\|^2+(k_n^2-1)M-2\langle v, k_n^2z_n-x_n\rangle)\},\\
Q_n=\{v\in C\hbox{\rm :}\ \langle x_0-x_n, x_n-v\rangle\geq 0\},\\
x_{n+1}=P_{C_n \cap Q_n}x_0,
\end{cases}
\end{equation*}
where  $M$ is an appropriate constant such that $M>\|v\|^2$ for
each $v\in C.$ Then $\{x_n\}$ converges to some $q= P_{F(T)}x_0.$
\end{theo}

\begin{proof}
Note that $T$ has a fixed point in $C$ \cite{10}. The key is to
show that if $T$ is asymptotically nonexpansive, then $T$ is also
relatively asymptotically nonexpansive. Take $p\in\hat{F}(T)$.
There exists a sequence $\{x_n\}\subset C$ such that
$x_n\rightharpoonup p$ and
$\lim_{n\rightarrow\infty}\|Tx_n-x_n\|$. Since $T$ is a
asymptotically nonexpansive, it is well-known that $T$ is
demiclosed, which yields $p\in F(T)$. On the other hand, we have
$F(T)\subset\hat{F}(T)$. In Hilbert spaces we know (1.5) reduces
to $\phi(x,y)=\|x-y\|^2, \ x,y\in H$. That is,
$\phi(T^nx,T^ny)\leq k_n^2\phi(x,y)$ is equivalent to
$\|T^nx-T^ny\|\leq k_n\|x-y\|$. Therefore, $T$ is also relatively
asymptotically nonexpansive. By using Theorem 2.1, it is easy to
obtain the desired conclusion. This completes the\break
proof.\hfill $\Box$
\end{proof}

\begin{theo}[\!]\hskip -.5pc {\bf \cite{19}.} \ \
Let $C$ be a closed convex subset of a Hilbert space H and let
$T\hbox{\rm :}\ C \rightarrow C$ be a nonexpansive mapping such
that $F(T)\neq \emptyset$. Assume that $\{\alpha_n\}_{n=0}^\infty$
and $\{\beta_n\}_{n=0}^\infty$ are sequences in $[0, 1]$ such that
$\alpha_n\leq 1-\delta$ for some $\delta\in (0, 1]$ and
$\beta_n\rightarrow 1$. Define a sequence $\{x_n\}_{n=0}^\infty$
in $C$ by the algorithm{\rm :}
\begin{equation*}
\begin{cases}
x_0\in C \ \ \ \mbox{arbitrarily\ chosen,}\\[.4pc]
z_n=\beta_nx_n+(1-\beta_n)Tx_n,\\[.4pc]
y_n =\alpha_nx_n+(1-\alpha_n)Tz_n, \\[.4pc]
C_n=\{v\in C\hbox{\rm :}\ \|y_n-v\|^2\leq\|x_n-v\|^2\\
\qquad\,\, +(1-\alpha_n)(\|z_n\|^2-\|x_n\|^2+2\langle x_n-z_n, v\rangle)\},\\[.4pc]
Q_n=\{v\in C\hbox{\rm :}\ \langle x_0-x_n, x_n-v\rangle\geq 0\},\\[.4pc]
x_{n+1}=P_{C_n \cap Q_n}x_0.
\end{cases}
\end{equation*}
Then $\{x_n\}$ converges in norm to $P_{F(T)}x_0.$
\end{theo}

\begin{proof}
It is well-known that the nonexpansive map is an asymptotically
nonexpansive map when $k_n=1$. By using Theorem 3.1, it is easy to
obtain the desired conclusion. This completes the proof.\hfill
$\Box$
\end{proof}

\begin{theo}[\!]\hskip -.5pc {\bf \cite{15}.} \ \
Let $C$ be a closed convex subset of a Hilbert space H and let
$T\hbox{\rm :}\ C \rightarrow C$ be a nonexpansive mapping such
that $F(T)\neq \emptyset$. Assume that $\{\alpha_n\}_{n=0}^\infty$
is a sequence in $(0, 1)$ such that $\alpha_n\leq 1-\delta$ for
some $\delta\in (0, 1]$. Define a sequence $\{x_n\}_{n=0}^\infty$
in $C$ by the algorithm{\rm :}
\begin{equation*}
\begin{cases}
x_0\in C \ \ \ \mbox{arbitrarily,}\\[.4pc]
y_n =\alpha_nx_n+(1-\alpha_n)Tx_n, \\[.4pc]
C_n=\{v\in C\hbox{\rm :}\ \|y_n-v\|\leq\|x_n-v\|,\\[.4pc]
Q_n=\{v\in C\hbox{\rm :}\ \langle x_0-x_n, x_n-v\rangle\geq 0\},\\[.4pc]
x_{n+1}=P_{C_n \cap Q_n}x_0.
\end{cases}
\end{equation*}
Then $\{x_n\}$ converges in norm to $P_{F(T)}x_0.$
\end{theo}

\begin{proof}
By taking $\beta_n=1$ in Theorem 3.2, we can obtain the desired
conclusion.\hfill $\Box$
\end{proof}

\end{document}